\documentclass[12 pt]{amsart}

\usepackage{amssymb,latexsym,amscd}
\usepackage[mathscr]{euscript}
\usepackage{pstricks}
\usepackage[all]{xy}
\usepackage{layout}

\title[Cohomology of Toroidal Orbifold Quotients]
{Cohomology of Toroidal Orbifold Quotients}

\newtheorem{theorem}{Theorem}

\newtheorem{lemma}[theorem]{Lemma}

\newtheorem{corollary}[theorem]{Corollary}
\newtheorem{definition}[theorem]{Definition}

\def\Proof{\medskip\noindent{\bf Proof: }}

\def\Z{\mathbb{Z}}

\def\Q{\mathbb{Q}}
\def\R{\mathbb{R}}
\def\F{\mathbb{F}}
\def\S{\mathbb{S}}

\def\bop{\displaystyle\mathop{\bigoplus}}
\def\colim{\displaystyle\mathop{\text{colim}}}

\def\Hom{\text{Hom}}

\newpsobject{grilla}{psgrid}{subgriddiv=1,griddots=10,gridlabels=6pt}

\begin{document}

\author[A.~Adem]{Alejandro Adem$^{*}$}
\address{Department of Mathematics,
University of British Columbia, Vancouver BC V6T 1Z2, Canada}
\email{adem@math.ubc.ca}
\thanks{$^{*}$Partially supported by NSERC}

\author[A.~N.~Duman]{Ali Nabi Duman}
\address{Department of Mathematics,
University of British Columbia, Vancouver BC V6T 1Z2, Canada}
\email{alid@math.ubc.ca}

\author[J.~M.~G\'omez]{Jos\'e Manuel G\'omez}
\address{Department of Mathematics,
University of British Columbia, Vancouver BC V6T 1Z2, Canada}
\email{josmago@math.ubc.ca}

\begin{abstract}
Let $\varphi:\Z/p\to GL_{n}(\Z)$ denote an integral representation
of the cyclic group of prime order $p$. This induces a $\Z/p$-action
on the torus $X=\R^{n}/\Z^{n}$. The goal of this paper is to
explicitly compute the cohomology groups $H^{*}(X/\Z/p;\Z)$ for any
such representation. As a consequence we obtain an explicit
calculation of the integral cohomology of the classifying space
associated to the family of finite subgroups for any
crystallographic group $\Gamma =\Z^n\rtimes\Z/p$  with prime
holonomy.

\end{abstract}

\maketitle

\section{Introduction}

Let $G$ be a finite group and $\varphi:G\to GL_{n}(\Z)$ an integral
representation of $G$. In this way $G$ acts linearly on $\R^{n}$
preserving the integral lattice $\Z^{n}$, hence inducing a
$G$-action on the torus $X_{\varphi}=X:=\R^{n}/\Z^{n}$. The quotient
$X\to X/G$ has the natural structure of a global quotient orbifold
and is an example of what is often called a toroidal orbifold. The
goal of this paper is to explicitly compute the cohomology groups
$H^{*}(X/G;\Z)$ for the particular case where $G=\Z/p$ is the cyclic
group of prime order.

\medskip

The indecomposable integral representations of $\Z/p$ have been 
completely classified (see \cite{CR}). In general, if $L$ is a 
$\Z G$-lattice then there are unique integers $r,s$ and $t$ and 
an isomorphism
\begin{equation*}
L\cong (\bop_{r} A_{i})\oplus (\bop_{s} P_{j})\oplus (\bop_{t}\Z),
\end{equation*}
where $\Z$ is the trivial $\Z/p$-module, each $A_{i}$ is an 
indecomposable module that corresponds to an element of the ideal 
class group of $\Z[\zeta_{p}]$, where $\zeta_{p}$ is a primitive 
$p$-th root of unity, and each $P_{j}$ is a projective indecomposable 
module of rank $p$ as an abelian group. In this case $L$ is said to be 
a $\Z/p$-module of type $(r,s,t)$. The situation simplifies after 
localizing at the prime $p$. Let $\Z_{(p)}$ be the ring of integers 
localized at the prime $p$. Then (see \cite{CR}), there are 
only three distinct isomorphism classes of indecomposable 
$\Z_{(p)}G$-lattices, namely the trivial module $\Z_{(p)}$, 
the augmentation ideal $IG$ and the group ring $\Z_{(p)} G$. 
Moreover, if $L$ is any finitely generated $\Z G$-lattice, then 
there is a $\Z G$-lattice 
$L'\cong IG^{r}\oplus \Z G^{s}\oplus  \Z^{t}$ and a 
$\Z G$-homomorphism $f:L'\to L$ such that $f$ is an isomorphism 
after tensoring with $\Z_{(p)}$.
\medskip

In this paper it is shown that given a $\Z G$-lattice $L$ induced 
by an integral representation $\varphi$, the cohomology groups 
$H^{*}(X_{\varphi}/G;\Z)$ only depends on the type of $L$. 
Moreover, if $L$ is of type $(r,s,t)$ then explicit descriptions 
for these cohomology groups are obtained in terms of $r,s$ and 
$t$. More precisely the goal of this paper is to prove the 
following theorem.

\begin{theorem}\label{maintheorem}
Let $G=\Z/p$ where $p$ is a prime number. Suppose that $X$ is 
the $G$-space induced by a $\Z G$-lattice $L$ of type $(r,s,t)$ 
and rank $n$. Then
\begin{displaymath}
H^{k}(X/G;\Z)\cong \Z^{\alpha_{k}}
\oplus (\Z/p)^{\beta_{k}},
\end{displaymath}
and the coefficients $\alpha_{k}$ and $\beta_{k}$ are given 
as follows: consider the formal power series in $x$
\[
F_{L}(x)=\left(\frac{1-(\alpha x)^{p}}{1-\alpha x}\right)^{r}
(1+\epsilon_{p} x^{p} )^{s}(1+x)^{t}
\]
subject to the relations $\alpha^{2}=1$, $\epsilon_{2}=\alpha$ 
and $\epsilon_{p}=1$ for $p>2$.  Using these relations, $F_{L}(x)$ 
can be written in the form
\[
F_{L}(x)=\sum_{i\ge 0}f_{i}x^{i}+\sum_{i\ge 0}g_{i}\alpha x^{i},
\]
then 
\[
\alpha_{k}=\frac{1}{p}\left[\binom{{n}}{{k}}+(p-1)(f_{k}-g_{k})
\right]. 
\]
Similarly, $\beta_{k}$ is obtained by writing the formal series in $x$ 
\[
T_{L}(x)=\frac{x(1+x)^{t}}{1-x^{2}}\left[
p^{r}x^{2}(1+x)^{s}-x^{2}+1-(1+\alpha x)(1+\epsilon_{p} x^{p} )^{s}
\left(\frac{1-(\alpha x)^{p}}{1-\alpha x}\right)^{r}\right]
\]
in the form $T_{L}(x)=\sum_{i\ge 0}\beta_{i}x^{i}+
\sum_{i\ge 0}\gamma_{i}\alpha x^{i}$, where $\alpha$ and 
$\epsilon_{p}$ subject to the same relations as above.
\end{theorem}
\medskip

Toroidal orbifolds naturally appear in different geometric 
contexts. In dimension six they are a source of key examples 
connected to mathematical aspects of orbifold string theory 
(see \cite{AGPP}). They also arise in the context of spaces 
of representations; it can be shown that the moduli space of 
isomorphism classes of flat connections on principal stable 
symplectic bundles over the torus $(\S^{1})^{n}$, can be 
described as the infinite symmetric product of a toroidal 
orbifold,
$$Rep(\Z^n, Sp):= \colim_{m\to \infty} 
Rep(\Z^n, Sp(m))\cong SP^{\infty}((\S^{1})^{n}/\Z/2)$$
where $\Z/2$ acts diagonally by complex conjugation, an action 
which arises from the direct sum of copies of the sign 
representation. In fact this space turns out to be a product 
of Eilenberg-MacLane spaces determined precisely by the 
homology of the quotient orbifold $(\S^{1})^{n}/\Z/2$ 
(see \cite{ACG} for details).

Similarly, recall that given a topological space $Y$, the $m$-th
cyclic product of $Y$ is defined to be the quotient
$CP^{m}(Y):=Y^{m}/\Z/m$, where $\Z/m$ acts on $Y^{m}$ by a 
cyclic permutation. The calculations here provide a complete 
computation for the homology of the $p$-th cyclic powers of 
any torus $(\S^1)^n$, as the permutation action corresponds to 
a direct sum of copies of the regular representation of $\Z/p$. 
Note that a method for such calculations was formulated long 
ago by Swan (\cite{Swan}); the approach outlined here is of 
course much more explicit.

\medskip

However the most important motivation for these calculations 
arises from the study of topological invariants of crystallographic 
groups with prime holonomy. Given a rank $n$ integral representation 
of $\Z/p$, it can be easily seen that this gives rise to an action 
of the semi-direct product $\Gamma=\Z^n\rtimes\Z/p$ on $Y=\R^n$ 
with the following crucial properties: for a subgroup 
$H\subset\Gamma$, $Y^H$ is non-empty if and only if $H$ is a 
finite subgroup of $\Gamma$, and furthermore in that case $Y^H$ 
is contractible. Thus $Y$ is a universal space for the family 
of finite subgroups in $\Gamma$, denoted by $\underline{E}\Gamma$ 
(see \cite{Luck} for definitions), and the associated classifying 
space is $\underline{B}\Gamma=\underline{E}\Gamma/\Gamma$, which in 
this case is precisely the orbifold quotient $(\S^1)^n/\Z/p$. Thus 
our main result together with the results in \cite{Charlap} 
can be reformulated as follows:

\begin{theorem}
Let $\Gamma$  be a crystallographic group with holonomy 
of prime order $p$, expressed as an extension
\begin{equation}\label{extension}
1\to L\to \Gamma \to \Z/p\to 1.
\end{equation}
Then the cohomology of the classifying space for the family 
of all finite subgroups in
$\Gamma$ can be explicitly computed and depends only on the
representation type of $L$ over the ring of integers localized at
$p$ as follows:
\begin{itemize}
\item If $\Gamma$ is torsion-free then $\Gamma$ is a Bieberbach 
group, $\underline{B}\Gamma= B\Gamma$ and 
$H^{*}(\underline{B}\Gamma;\Z)$ can be computed using 
\cite[Theorem 2]{Charlap}. 

\item If $\Gamma$ is not torsion-free, then the sequence 
(\ref{extension}) splits, $\Gamma=L\rtimes \Z/p$ and 
$H^{*}(\underline{B}\Gamma;\Z)$ can be computed using 
Theorem \ref{maintheorem}.
\end{itemize}
\end{theorem}

Using methods from equivariant $K$--theory and an analysis analogous
to that done in this paper, these computations will serve as
important input for the calculation of the complex $K$--theory of
$B\Gamma$ for $\Gamma$ a crystallographic group with prime holonomy.

\section{Preliminaries}

Let $G$ be a finite group and $\varphi:G\to GL_{n}(\Z)$ an integral
representation of $G$. Consider $X=X_{\varphi}$ the $G$-space 
induced by the representation $\varphi$. Then the fibration sequence
\begin{equation}\label{fibration sequence}
X\to X\times_{G} EG\to BG
\end{equation}
induces a long exact sequence of homotopy groups
\[
\cdots \to\pi_{i}(X)\to \pi_{i}(X\times_{G} EG)\to \pi_{i}(BG)\to
\pi_{i-1}(X)\to\cdots.
\]
This sequence is trivial for $i>1$ and thus $X\times_{G} EG$ is an
Eilenberg-Maclane space of type $K(\Gamma,1)$, where
$\Gamma:=\pi_{1}(X\times_{G} EG)$ fits into a short exact sequence
\begin{equation}\label{seqhomotopygroups}
1\to \pi_{1}(X)\to \Gamma \to G\to 1.
\end{equation}
The action of $G$ on $X$ makes
\[
L:=\pi_{1}(X)\cong H_{1}(X;\Z)\cong \Z^{n}
\]
into a $\Z G$-module that corresponds to the representation 
$\varphi$. Moreover, $[0]\in \R^{n}/\Z^{n}=X$ is a fixed point 
for this action and thus (\ref{fibration sequence}) has a section. 
This implies that the extension (\ref{seqhomotopygroups}) splits 
and thus $\Gamma\cong L\rtimes G$. For example, when the 
representation $\varphi$ is faithful the group $\Gamma$ is a 
crystallographic group.

\medskip

The cohomology groups of groups of the form
$\Gamma\cong L\rtimes G$, for $G=\Z/p$ with $p$ a prime number, 
were computed in \cite[Theorem 1.1]{AGPP}. In there it was proved 
that the Lyndon-Hochschild-Serre spectral sequence associated 
to the short exact sequence
\[
1\to L\to \Gamma\to G\to 1
\]
collapses on the $E_{2}$-term without extension problems. This can 
be seen as follows. Suppose first that $L$ is a $\Z G$-lattice 
of the form $L=IG^{r}\oplus \Z G^{s}\oplus  \Z^{t}$. For such 
lattices it follows by \cite[Theorem 3.2]{AGPP} and 
\cite[Proposition 3.3]{AGPP} that there is a special free 
resolution $\epsilon:F\to \Z$ of $\Z$ as a $\Z[L]$-module 
admitting an action of $G$ compatible with $\varphi$. 
Thus by \cite[Theorem 2.4]{AGPP} the corresponding 
Lyndon-Hochschild-Serre spectral sequence collapses in 
this particular case. Suppose now that $L$ is any 
$\Z G$-lattice. Then we can find a $\Z G$-lattice 
$L'\cong IG^{r}\oplus \Z G^{s}\oplus  \Z^{t}$ and a 
$\Z G$-homomorphism $f:L'\to L$ such that $f$ is an isomorphism 
after tensoring with $\Z_{(p)}$. By comparing the spectral sequences 
corresponding to $L$ and $L'$ as done in \cite[Theorem 4.1]{AGPP} 
it can be seen that the Lyndon-Hochschild-Serre corresponding 
to $L$ also collapses on the $E_{2}$-term without extension 
problems. Therefore, for any $k\ge 0$
\[
H^{k}(\Gamma;\Z)\cong \bigoplus_{i+j=k}
H^{i}(G;\wedge^{j}(L^{*})).
\]
Here $L^{*}$, as usual, denotes the dual $G$-module $\Hom(L,\Z)$.  
As an application of this, by \cite[Theorem 1.2]{AGPP}, if 
$G=\Z/p$ acts on $X$ via a representation
$\varphi:G\to GL_{n}(\Z)$, then for each $k\ge 0$
\begin{equation}\label{collapseSerre}
H^{k}_{G}(X;\Z)\cong \bigoplus_{i+j=k}H^{i}(G;H^{j}(X;\Z)).
\end{equation}

This completely describes the additive structure of the 
equivariant cohomology groups $H^{k}_{G}(X;\Z)$. Moreover, these 
groups can explicitly be computed as observed in \cite{Charlap} 
in the following way. Suppose that $L$ is a $\Z G$-lattice of 
type $(r,s,t)$ and rank $n$. Consider the formal power series 
in $x$
\[
F_{L}(x)=\left(\frac{1-(\alpha x)^{p}}{1-\alpha x}\right)^{r}
(1+\epsilon_{p} x^{p} )^{s}(1+x)^{t}
\]
subject to the relations $\alpha^{2}=1$, $\epsilon_{2}=\alpha$ 
and $\epsilon_{p}=1$ for $p>2$.  Using these relations, the 
formal series $F_{L}(x)$ can be written in the form
\[
F_{L}(x)=\sum_{i\ge 0}f_{i}x^{i}+\sum_{i\ge 0}g_{i}\alpha x^{i}
\]
for integer numbers $f_{i}$ and $g_{i}$ for $i\ge 0$. The 
coefficients $f_{i}$ and $g_{i}$ determine the type of the 
$\Z G$-module $\wedge^{i} L$. Indeed, by 
\cite[Corollary 5.8]{Charlap} if  
\[
h_{i}=\frac{1}{p}\left[\binom{{n}}{{i}}+(p-1)(f_{i}-g_{i})\right] 
\]
then $\wedge^{i} L$ is of type $(g_{i},h_{i}-f_{i},f_{i})$ 
and thus
\[
H^{0}(G,\wedge^{i}L)=\Z^{h_{i}}, \ \
H^{1}(G,\wedge^{i}L)=(\Z/p)^{g_{i}}, \ \
H^{2}(G,\wedge^{i}L)=(\Z/p)^{f_{i}}.
\]

As a corollary the following is obtained.

\begin{corollary}\label{equivariant cohomology}
Suppose $G=\Z/p$ acts on $X$ via a representation
$\varphi:G\to GL_{n}(\Z)$ inducing a $\Z G$-lattice $L$ of 
type $(r,s,t)$. Then for each $k\ge 0$
\[
H^{k}_{G}(X;\Z)\cong \Z^{a_{k}}\oplus (\Z/p)^{b_{k}}
\]
and the coefficients $a_{k}$ and $b_{k}$ are given as follows:
write 
\[
F_{L}(x)=\sum_{i\ge 0}f_{i}x^{i}+\sum_{i\ge 0}g_{i}\alpha x^{i}
\]
then
\[
a_{k}=\frac{1}{p}\left[\binom{{n}}{{k}}+(p-1)(f_{k}-g_{k})\right]. 
\]
Also, $b_{k}$ is obtained by writing the formal series in $x$ 
\[
G_{L}(x)=\alpha x \left(\frac{1-(\alpha x)^{k}}{1-\alpha x}\right)
\left(\frac{1-(\alpha x)^{p}}{1-\alpha x}\right)^{r}
(1+\epsilon_{p} x^{p} )^{s}(1+x)^{t}
\]
in the form $G_{L}(x)=\sum_{i\ge 0}b_{i}x^{i}+
\sum_{i\ge 0}c_{i}\alpha x^{i}$ where $\alpha$ and $\epsilon_{p}$ 
subject to the same relations 
as above.
\end{corollary}

This corollary can be used as a first step towards the 
computation of the cohomology groups of the form $H^{k}(X/G;\Z)$ 
as is shown next. 

\begin{corollary}\label{free part}
Let $X$ be induced by an integral representation 
$\varphi:G\to GL_{n}(\Z)$ inducing a $\Z G$-lattice of type 
$(r,s,t)$. Then 
\[
H^{k}(X/G;\Z)=F_{k}\oplus T_{k},
\]
where $F_{k}$ is a free abelian group of rank 
\[
\alpha_{k}=\frac{1}{p}\left[\binom{{n}}{{k}}+(p-1)(f_{k}-g_{k})
\right]
\]
and $T_{k}$ is a $p$-torsion abelian group.
\end{corollary}
\Proof
Consider the map $\phi:X\times_{G}EG \to X/G$ obtained by  
mapping $EG$ to a point. Given $x\in X$, $\phi^{-1}([x])\cong BG_{x}$, 
where $G_{x}$ is the isotropy group of $G$ at $x$. Since $G=\Z/p$, 
then $G_{x}$ is either the 
trivial subgroup or $G$. In either case $BG_{x}$ has trivial 
cohomology with coefficients in $\Q$ and also with coefficients in 
$\F_{q}$, the field with $q$-elements for a prime $q$ different from 
$p$. The Vietoris-Begle Theorem shows that $\phi$ induces isomorphisms
\begin{align*}
\phi^{*}:H^{*}(X/G;\Q)&\stackrel{\cong}{\rightarrow} 
H^{*}_{G}(X;\Q),\\
\phi^{*}:H^{*}(X/G;\F_{q})&\stackrel{\cong}{\rightarrow} 
H^{*}_{G}(X;\F_{q}).
\end{align*}
Using the previous corollary and the universal 
coefficient theorem we see that 
\[
H^{k}_{G}(X;\Q)\cong \Q^{\alpha_{k}} \text{ and } 
H^{k}_{G}(X;\F_{q})\cong \F_{q}^{\alpha_{k}}
\]
and the corollary follows.
\qed
\medskip

The previous corollary reduces our problem to computing 
$T_{k}$, the $p$-torsion subgroup of $H^{k}(X/G;\Z)$. The strategy 
that we will follow to  compute these groups is as follows. Consider
$\phi:X\times_{G}EG\to X/G$ the map defined above and 
let $F$ be the fixed point set of $X$ under 
the action of $G$. Then \cite[Proposition VII 1.1]{Bredon} shows 
that $\phi$ induces an isomorphism
\[
\phi^{*}:H^{*}(X/G,F;\Z)\to H^{*}_{G}(X,F;\Z).
\]
Via this isomorphism, the groups $T_{k}$ can be computed using  
the following steps. First we compute the $p$-torsion subgroups of 
$H^{*}_{G}(X,F;\Z)$. Then we use this information together 
with the long exact sequence in cohomology 
associated to the pair $(X/G,F)$ to deduce the structure of $T_{k}$. 
\newline

We establish now some notation. Let $R^{*}=H^{*}(G,\Z)$. 
Then $R^{*}$ can be seen as a graded 
commutative ring whose structure is given by $R^{*}=\Z[t]/(pt)$, 
where $deg(t)=2$. Graded $R^{*}$-modules of the form 
$M=\bigoplus_{n\ge 0} M^{n}$, where $M^{n}$ is a finite 
dimensional $\F_{p}$-vector space for $n> 0$, appear naturally in 
our computations. For such modules we have the following definition.

\begin{definition}\label{definition power series}
Given a graded $R^{*}$-module $M=\bigoplus_{n\ge 0} M^{n}$ as above, 
define the formal power series in $\Z[[x]]$
\[
q_{M}(x):=a_{1}x+a_{2}x^{2}+a_{3}x^{3}+\cdots 
\]
where $a_{n}=\text{dim}_{\F_{p}}M^{n}$ for $n>0$. 
\end{definition}

If $M$ is a graded $R^{*}$-module $M$ as above, the series 
$q_{M}(x)$ together with $M^{0}$ completely determine 
the structure of $M$ as an abelian group. For example, 
if $L$ is a $\Z G$-lattice, then $H^{*}(G,L)$ is a graded 
$R^{*}$-module of this kind and the series $q_{H^{*}(G,L)}(x)$ 
can be explicitly computed in the following way. Suppose first 
that $A$ is an indecomposable module of rank $p-1$ that 
corresponds to an element of the ideal class group of 
$\Z[\zeta_{p}]$. Then by \cite[Corollary 1.7]{Adem}, it follows 
that $S^{*}:=H^{*}(G,A)$ is a graded $R^{*}$-module such that 
$S^{n}=0$ for $n$-even and $S^{n}=\Z/p$ for $n>0$ odd. Also, 
given any projective indecomposable module $P$ of rank $p$ 
then by \cite[Proposition 1.5]{Adem}, $H^{*}(G,P)=\Z$ is the 
trivial graded $R^{*}$-module concentrated on degree $0$. 
Therefore, given a $\Z G$-lattice $L$ of type 
$(r,s,t)$ there is an isomorphism of graded 
$R^{*}$-modules
\begin{equation*}
H^{*}(G,L)\cong \Z^{{s}}\oplus(R^{*})^{t}\oplus (S^{*})^{r}.
\end{equation*}
In particular 
\begin{align}\label{Series for (r,s,t)}
q_{H^{*}(G,L)}(x)&=rx+tx^{2}+rx^{3}+tx^{4}\cdots\\
&=\frac{rx+tx^{2}}{1-x^{2}}.
\end{align}

To compute the $p$-torsion subgroup of $H^{k}_{G}(X,F;\Z)$ 
we use the Serre spectral sequence
\begin{equation}\label{SerreSS}
E^{i,j}_{2}=H^{i}(G,H^{j}(X,F;\Z))\Longrightarrow 
H^{i+j}_{G}(X,F;\Z)
\end{equation}
associated to the pair $(X,F)$.  In this spectral sequence when  
$j$ and $r$ are fixed, $E_{r}^{*,j}$ is a graded $R^{*}$-module of 
the kind considered in Definition \ref{definition power series} and 
each differential
\[
d_{r}:E_{r}^{*,j}\to E_{r}^{*,j-r+1}
\]
is homomorphism of graded $R^{*}$-modules of degree $r$. In here 
we will show that  the different formal power series 
$q_{E^{*,j}_{\infty}}(x)$, for $j\ge 0$, determine the 
$p$-torsion subgroups of $H^{k}_{G}(X,F;\Z)$. 
This will be done by determining the nontrivial 
differentials in the spectral sequence (\ref{SerreSS}). To do this  
we first consider the particular cases of $\Z G$-lattices of 
type $(r,0,0)$ in section \ref{(r,0,0)} and type $(0,s,0)$ in section 
\ref{(0,s,0)}. Then we use this information to handle the 
general case in section \ref{general case}.

\medskip

The following lemma plays a key role in our computations.

\begin{lemma}\label{lemmaBredon}
Suppose that $p$ is a prime number. Let $G=\Z/p$ act on a finistic 
space $X$ with fixed point set $F$. If there is an integer $N$ such 
that $H^{k}(X,F;\Z)=0$ for $k>N$, then $H^{k}_{G}(X,F;\Z)=0$ for 
$k>N$.
\end{lemma}
\Proof
This follows by applying \cite[Exercise III.9]{Bredon} and
\cite[Proposition VII 1.1]{Bredon}.
\qed
\medskip

\section{Structure of the fixed points}

In this section we investigate the nature of the fixed point set 
of the action of $G$ on a torus $X$ induced by a general 
$\Z G$-lattice $L$. 

\medskip

To start, note that if $L$ and $M$ are two $\Z G$-lattices then
as $G$-spaces
\begin{equation}\label{additiveprop}
X_{L\oplus M}=X_{L}\times X_{M}.
\end{equation}
In particular this shows that $F_{L\oplus M}=F_{L}\times F_{M}$. 
Here $F_{L}$ denotes the fixed point set of the action of $G$ 
on $X_{L}$ for a given $\Z G$ lattice $L$. 

We consider next the particular cases $L=A$, $L=P$ and $L=\Z$, 
where $A$ and $P$ are indecomposable modules as of rank $p-1$ 
and rank $p$ respectively, as described before. 

\begin{lemma}\label{fixed point set case IG}
Let $A$ be an indecomposable module of rank $p-1$ corresponding 
to an element of the ideal class group of $\Z[\zeta_{p}]$. Then the 
fixed point set $F$ of the $G$-action on the induced torus $X$ is a 
discrete set with $p$ points.
\end{lemma}
\Proof
Consider the short exact sequence of $G$-modules defining
the $G$-space $X$
\[
0\to A\to A\otimes \R\to (A\otimes \R)/A=X\to 0.
\]
This short exact sequence induces a long exact sequence on the level
of group cohomology
\[
0\to H^{0}(G,A)\to H^{0}(G,A\otimes \R)\to H^{0}(G,X)\to
H^{1}(G,A)\to H^{1}(G,A\otimes \R)\to \cdots.
\]
Note that $H^{1}(G,A\otimes \R)=0$ and
$H^{0}(G,A)= H^{0}(G,A\otimes \R)=0$, thus there is an isomorphism
\[
F=H^{0}(G,X)\cong H^{1}(G,A)\cong \Z/p.
\]
\qed
\medskip

\begin{lemma}\label{cohopairsZG} 
Let $P$ be a projective indecomposable 
module of rank $p$. Then $F_{P}\cong \S^{1}$. Moreover, 
there is a commutative diagram 
\begin{equation}\label{fixpoints(0,s,0)}
\xymatrix{
\S^{1}\ar[d]_{h_{1}}\ar[r]^{\Delta}&(\S^{1})^{p}\ar[d]^{h_{2}}\\
F_{P}\ar[r]^{i}&X_{P}.
}
\end{equation} 
where $\Delta$ denotes the diagonal inclusion of  $ \S^{1}$ into 
$(\S^{1})^{p}$ and $h_{1}$, $h_{2}$ are covering maps of degree 
relatively prime with $p$.
\end{lemma}
\Proof
Consider the exact sequence defining $X_{P}$
\[
0\to P\to P\otimes \R \to X_{P}\to 0.
\]
This yields a short exact sequence of the form
\[
0\to P^{G}\to (P\otimes \R)^{G} \to (X_{P})^{G}\to 0
\]
as $H^{1}(G,P)=0$. Note that $P^{G}\cong \Z$ and 
$(P\otimes \R)^{G}\cong \R$, therefore $F=X_{P}^{G}\cong \S^{1}$. 
To prove the second assertion, suppose that $L=\Z G$. Then in this 
particular case it is easy to see that the inclusion 
$i:F_{\Z G}\to X_{\Z G}$ corresponds to the the diagonal 
inclusion $\Delta:\S^{1}\to (\S^{1})^{p}$. Let 
$f:\Z G\to P$  be homomorphism of $\Z G$-modules such that
$f$ is an isomorphism after tensoring with $\Z_{(p)}$. In particular, 
$f$ is an isomorphism after tensoring with $\R$ and we have 
a commutative diagram 
\[
\xymatrix{
0\ar[r] &\Z\ar[d]_{f^{G}}\ar[r]& \R\ar[r]
\ar[d]_{(f\otimes 1)^{G}}& (X_{\Z G})^{G}
\ar[r]\ar[d]_{\tilde{f}} &0\\
0\ar[r]  &P^{G}\ar[r]& (P\otimes \R)^{G}\ar[r]& (X_{P})^{G}\ar[r] &0.
}
\]

The map $f^{G}:\Z=(\Z G)^{G}\to P^{G}\cong \Z$ must be 
multiplication by a number $q$, where $q$ is relatively prime to $p$. 
This proves that the map 
\[
h_{1}:\S^{1}\to F_{P}
\]
induced by $f$ is a degree $q$ covering map. Similarly, 
$h_{2}:(\S^{1})^{p}\to X_{P}$ is a covering map of degree 
relatively prime to $p$.
\qed
\medskip

Notice that when $L=\Z$ is the trivial $G$-module then 
$X_{\Z}=F_{\Z}=\S^{1}$. As a corollary we obtain.

\begin{corollary}\label{fixedpointset}
Let $X$ be the $G$-space induced by a $\Z G$-lattice of type 
$(r,s,t)$. Then 
\[
F:=X^{G}\cong \bigsqcup_{p^{r}}(\S^{1})^{s+t}.
\]
\end{corollary} 

\begin{lemma}\label{isocohomologypairs}
Let $L$ and $M$ be two $\Z G$-lattices. Assume that $f:L\to M$ is 
a $\Z G$-homomorphism that is an isomorphism after tensoring with 
$\Z_{(p)}$. Then $f$ induces an isomorphism
\[
f^{*}:H^{*}_{G}(X_{M},F_{M};\Z_{(p)})\to H^{*}_{G}(X_{L},F_{L};
\Z_{(p)}).
\]
\end{lemma}
\Proof 
Let $K=\Z_{(p)}$. Consider the Serre spectral sequence with 
$K$-coefficients
\[
E_{2}^{i,j}=H^{i}(G,H^{j}(X_{M},F_{M};K))\Longrightarrow 
H^{i+j}_{G}(X_{M},F_{M};K). 
\]
Similarly we obtain a spectral sequence $\tilde{E}_{r}^{i,j}$ 
associated to $L$. The map $f:L\to M$ induces a map of spectral 
sequences
\[
f_{r}^{i,j}:E_{r}^{i,j}\to \tilde{E}_{r}^{i,j}.
\] 
We will prove the lemma by showing that $f$ induces an isomorphism 
on the corresponding $E_{2}$-terms. Let $L_{K}^{*}=L^{*}\otimes K$ 
and $M_{K}^{*}=M^{*}\otimes K$. To start note that 
$H^{j}(X_{L};K)\cong \wedge^{j}L_{K}^{*}$ and 
similarly $H^{j}(X_{M};K)\cong \wedge^{j}M_{K}^{*}$. By 
hypothesis $\wedge^{j}L_{M}^{*}$ and $\wedge^{j}L_{M}^{*}$ are 
isomorphic $KG$--modules, with isomorphism induced by $f$. Thus 
\begin{equation}\label{auxrelative1}
f^{*}:H^{j}(X_{M};K)\to H^{j}(X_{L};K)
\end{equation} 
is an isomorphism of $KG$-modules. On the other hand, since 
$f:L\to M$ is an isomorphism after tensoring with $K$, then $L$ 
can be seen as a sub-lattice of $M$ of finite index $q$, with $q$ 
relatively prime with $p$. Consider the map $f^{G}:F_{L}\to F_{M}$
induced by $f$ on the level of fixed points. An argument similar 
to that in Lemma \ref{cohopairsZG} can be used to show that 
$f^{G}$ is a covering map of degree relatively prime to $p$. 
In particular 
\begin{equation}\label{auxrelative2}
(f^{G})^{*}:H^{j}(F_{M};K)\to H^{j}(F_{L};K)
\end{equation}
is an isomorphism of $KG$-modules. Finally, note that $f$ induces 
a morphism between the long exact sequences in cohomology with 
$K$-coefficients associated to the pairs $(X_{L},F_{L})$ and 
$(X_{M},F_{M})$. By (\ref{auxrelative1}), (\ref{auxrelative2}) 
and the $5$-lemma we conclude that  
\[
f^{*}:H^{*}(X_{M},F_{M};K)\to H^{*}(X_{L},F_{L};K).
\]
is an isomorphism of $KG$-modules. This proves the lemma.
\qed
\medskip

\section{Modules of type $(r,0,0)$.}\label{(r,0,0)}

In this section we consider the particular case of a 
$\Z G$-lattice $L$ of type $(r,0,0)$ determined by an integral 
representation $\varphi$. 

\medskip

Suppose that $L$ is such a lattice. Then as an abelian group $L$ 
has rank $n=r(p-1)$ and thus the associated torus $X$ has 
rank $r(p-1)$. Consider the Serre spectral associated 
to the pair $(X,F)$
\begin{equation}\label{Serre1}
E_{2}^{i,j}=H^{i}(G,H^{j}(X,F;\Z))\Longrightarrow 
H_{G}^{i+j}(X,F;\Z).
\end{equation}
In this case $F$ is a finite set with $p^{r}$ points 
by Lemma \ref{fixedpointset}. In particular $H^{k}(X,F)=0$ 
if $k>r(p-1)$, thus  Lemma \ref{lemmaBredon} implies that 
$H^{k}_{G}(X,F)=0$ for $k>r(p-1)$. We will use this fact 
to show that in the spectral sequence (\ref{Serre1}) 
all the differentials are trivial except for the differentials 
of the form
\[
d_{j}:E_{j}^{*,j}\to E_{j}^{*,1},
\] 
whenever $1\le j\le r(p-1)$. 
Moreover, we will see that for such $j$ the differential 
$d_{j}$ is injective on positive degrees. To show this, note 
that as $F$ is a discrete set then $H^{0}(X,F)=0$ and 
for $j\ge 2$
\[
H^{j}(X,F)\cong H^{j}(X)\cong \wedge^{j}L^{*}.
\]
Since $L^{*}$ is of type $(r,0,0)$, if
\[
F_{L}(x)=\left(\frac{1-(\alpha x)^{p}}{1-\alpha x}\right)^{r}
=\sum_{i\ge 0}f_{i}x^{i}+\sum_{i\ge 0}g_{i}\alpha x^{i}
\]
and $j\ge 2$, then $H^{j}(X,F)$ is a $\Z G$-lattice of type 
$(g_{j},h_{j}-f_{j},f_{j})$ with 
\[
h_{j}=\frac{1}{p}\left[\binom{{n}}{{j}}+(p-1)(f_{j}-g_{j})
\right]. 
\]
In particular there is an 
isomorphism of graded $R^{*}$-modules
\begin{equation}\label{gradedmod1}
H^{*}(G,H^{j}(X,F))\cong\Z^{h_{j}-f_{j}}\oplus 
(R^{*})^{f_{j}}\oplus (S^{*})^{g_{j}}.
\end{equation}
Define $p_{r}(j)$ to be the number of all
possible sequences of integers $l_{1},...,l_{r}$ such that
$0\le l_{i}\le p-1$ and $l_{1}+\cdots +l_{r}=j$. Then in this 
case it is easy to see that for $j$ even $f_{j}=p_{r}(j)$ and 
$g_{j}=0$ and for $j$ odd $f_{j}=0$ and $g_{j}=p_{r}(j)$. 
Let's compute now $H^{*}(G,H^{1}(X,F))$. The long exact sequence in 
cohomology associated to the pair $(X,F)$ gives a short exact sequence 
of $\Z G$-modules
\[
0\to \Z^{p^{r}-1}\to H^{1}(X,F)\to H^{1}(X)\to 0.
\]
In particular, as a group, $H^{1}(X,F)$ is a free abelian group 
and there is a long exact sequence
\begin{equation}\label{ses1}
\cdots\to H^{i}(G, \Z^{p^{r}-1})\to H^{i}(G,H^{1}(X,F))\to 
H^{i}(G,H^{1}(X)) \stackrel{\delta}{\rightarrow} 
H^{i+1}(G, \Z^{p^{r}-1})\to\cdots.
\end{equation}
We claim that when $i$ is odd $H^{i}(G,H^{1}(X,F))=0$.
To see this it is enough to show that if $i$ is odd
and sufficiently large $H^{i}(G,H^{1}(X,F))=0$. Pick $i$ odd 
with $i>r(p-1)$ so that $H^{i+1}_{G}(X,F)=0$ for such $i$.
In particular there are no nontrivial permanent cocycles
in total degree $i+1$ in the spectral sequence (\ref{Serre1}).
Trivially all the differentials with source 
$E_{2}^{i,1}=H^{i}(G,H^{1}(X,F))$ are zero and
therefore any element in $E_{2}^{i,1}$ must be in the image of
some differential. However, (\ref{gradedmod1}) 
implies that any differential with target $E_{k}^{i,1}$ has a 
trivial source. This shows that $H^{i}(G,H^{1}(X,F))=0$ 
for $i$ odd. In particular, the long exact sequence (\ref{ses1})
reduces to the short exact sequence
\begin{equation}\label{sesgc}
0\to H^{i}(G,H^{1}(X))\to H^{i+1}(G,\Z^{p^{r}-1})\to
H^{i+1}(G,H^{1}(X,F))\to 0,
\end{equation}
for $i>0$ odd. Since $H^{1}(X)\cong L^{*}$ is a $\Z G$-module 
of type $(r,0,0)$ it follows that $H^{1}(X,F)$ is a $\Z G$-module 
of type $(0,r,p^{r}-r-1)$ and 
\begin{equation*}
H^{*}(G,H^{1}(X,F))\cong\Z^{r}\oplus (\bigoplus_{p^{r}-r-1}R^{*}).
\end{equation*}
This describes the $E_{2}$-term of the spectral sequence 
(\ref{Serre1}). Consider now the Serre spectral sequence
\[
\tilde{E}_{2}^{i,j}=H^{i}(G,H^{j}(X))\Longrightarrow 
H^{i+j}_{G}(X)
\]
associated to the fibration sequence $X\to X\times_{G} EG
\to BG$. As it was pointed out before this sequence collapses 
on the $E_{2}$-term. The inclusion $f:X\to (X,F)$ defines a map 
of spectral sequences $f^{i,j}_{k}:E_{k}^{i,j}\to \ 
\tilde{E}_{k}^{i,j}$. Notice that $f^{i,j}_{2}$ is an isomorphism 
when $j\ge 2$. This shows that the only possibly nontrivial 
differentials in (\ref{Serre1}) are of the form
\[
d_{j}:E_{j}^{*,j}\to E_{j}^{*,1}
\]
for $2\le j\le r(p-1)$. We determine next the nature of these 
differentials. Note that the factor $\Z^{h_{j}-f_{j}}\subset 
(H^{j}(X,F))^{G}$ lies in the image of the norm map. Consider the 
transfer map associated to the trivial subgroup $\{1\}
\hookrightarrow G$. This map preserves the filtrations that induce 
the Serre spectral sequence and thus it induces a map of the 
corresponding spectral sequences
\[
\tau_{\{1\}}^{G}:H^{i}(\{1\},H^{j}(X,F))\to 
H^{i}(G,H^{j}(X,F)).
\]
Since the image of the transfer map
$\tau_{\{1\}}^{G}:H^{0}(\{1\},H^{j}(X,F))\to H^{0}(G,H^{j}(X,F))$ 
consists of elements in the image of the norm map, it follows that 
all the differentials in the Serre spectral sequence 
(\ref{Serre1}) are trivial on the summand $\Z^{h_{j}-f_{j}}$. 
Let's show now that the differential $d_{j}$ is injective on positive 
degrees. To see this it suffices to show that $d_{j}:E_{j}^{i,j}\to 
E_{j}^{i+j,1}$ is injective for $i$ big enough. If $i>r(p-1)$ then 
$H_{G}^{i}(X,F)=0$, therefore in total degree $i$ with $i>r(p-1)$ there 
are no nontrivial permanent cocycles. Since all the differentials
landing in $E_{j}^{*,j}$ are trivial for $j\ge 2$ this forces
$d_{j}:E_{j}^{i,j}\to E_{j}^{i+j,1}$ to be injective when $i>r(p-1)$.

The above can be summarized in the following way. If $j\ne 1$ then 
\begin{equation}\label{Caser1}
q_{E_{\infty}^{*,j}}(x)=0.
\end{equation}
On the other hand, since $H^{1}(X,F)$ is of type $(0,r,p^{r}-r-1)$
then by (\ref{Series for (r,s,t)})
\[
q_{E_{2}^{*,1}}(x)=\frac{(p^{r}-r-1)x^{2}}{1-x^{2}}.
\]
Also, $d_{j}:E_{j}^{*,j}\to E_{j}^{*,1}$ is a 
homomorphism of graded $R^{*}$-modules of degree $j$ that is 
injective on positive degrees and $E_{j}^{*,j}=E_{2}^{*,j}$ is a 
$\Z G$-lattice of type $(g_{j},h_{j}-f_{j},f_{j})$, it follows that 
\[
q_{E_{j+1}^{*,1}}(x)=q_{E_{j}^{*,1}}(x)-
\left(\frac{f_{j}+xg_{j}}{1-x^{2}}\right)x^{j}.
\] 
Therefore
\begin{align}
q_{E_{\infty}^{*,1}}(x)&=\frac{(p^{r}-r-1)x^{2}}{1-x^{2}}-
\sum_{j\ge 2}\left(\frac{f_{j}+xg_{j}}{1-x^{2}}\right)x^{j}
\\&=\frac{1}{1-x^{2}}\left[(p^{r}-1)x^{2}+1
-\sum_{j\ge 0}(f_{j}+xg_{j})x^{j}\label{Caser2}
\right].
\end{align}

This completely characterizes the $E_{\infty}$-term in the 
spectral sequence (\ref{Serre1}) and determines $H_{G}^{*}(X,F;\Z)$ 
up to extension problems. Let's assume for a moment that there are no 
extension problems in this case. If this is true then 
\[
H^{k}_{G}(X,F;\Z)\cong \Z^{\kappa_{k}}\oplus (\Z/p)^{\lambda_{k}},
\] 
for some integers $\kappa_{k}$ and $\lambda_{k}$. Moreover, the 
integers $\lambda_{k}$ are determined by the formal series 
\[
\bar{Q}_{L}(x):=\sum_{j\ge 0}\lambda_{j}x^{j}
=\sum_{j\ge 0}x^{j}q_{E^{*,j}_{\infty}}(x).
\]
From (\ref{Caser1}) and (\ref{Caser2}) it follows that 
\[
\bar{Q}_{L}(x)=\sum_{j\ge 1}x^{j}q_{E_{\infty}^{*,j}}(x)
=\frac{x}{1-x^{2}}\left[p^{r}x^{2}-x^{2}+1
-\sum_{j\ge 0}(f_{j}+xg_{j})x^{j}
\right].
\]

Let's show now that indeed there are no extension problems in the 
spectral sequence (\ref{Serre1}). To see this we need 
to consider the Serre spectral sequence computing 
$H_{G}^{*}(X,F;\F)$, where  $\F=\Q$, $\F=\F_{p}$ and $\F=\F_{q}$, 
for a prime number $q$ different from $p$. Arguments 
similar to those formulated above can be applied in these 
cases to obtain explicit descriptions of the $E_{\infty}$-term 
in these spectral sequences. This way we obtain 
\begin{align*}
H_{G}^{k}(X,F;\Q)&\cong \Q^{\kappa_{k}},\\ 
H_{G}^{k}(X,F;\F_{q})&\cong \F_{q}^{\kappa_{k}},\\
H_{G}^{k}(X,F;\F_{p})&\cong 
\F_{p}^{\kappa_{k}+\lambda_{k}+\lambda_{k+1}}.
\end{align*}

The only way that this is possible is that indeed there 
are no extension problems in the spectral sequence (\ref{Serre1}).
To finish note that in Lemma \ref{long exact torsion} it 
is proved that the $p$-torsion of $H^{k}(X/G;\Z)$ is a 
$\F_{p}$-vector space of dimension $\beta_{k}$ that in this 
case equals $\lambda_{k}$. Thus these groups are determined 
by the formal power series
\[
Q_{L}(x):=\sum_{j\ge 0}\beta_{j}x^{j}=\bar{Q}_{L}(x).
\] 
To find the explicit description for the series $Q_{L}(x)$ 
notice that  
\[
(1+\alpha x)F_{L}(x)=(1+\alpha x)\left(\frac{1-(\alpha x)^{p}}
{1-\alpha x}\right)^{r}=\sum _{i\ge 0}(f_{i}+xg_{i})x^{i}+
\sum _{i\ge 0}\alpha(xf_{i}+g_{i})x^{i}.
\]
Thus if
\[
T_{L}(x)=\frac{x}{1-x^{2}}\left[
p^{r}x^{2}-x^{2}+1-(1+\alpha x)
\left(\frac{1-(\alpha x)^{p}}{1-\alpha x}\right)^{r}\right]
\]
then $T_{L}(x)=Q_{L}(x)+\alpha R_{L}(x)$, for some 
$R_{L}(x)\in \Z[[x]]$. This together with 
Corollary \ref{free part} prove the following theorem.

\begin{theorem}\label{theoremtype(r,0,0)}\footnote{This result 
forms part of the Doctoral dissertation of the second author 
(see \cite{ThesisAli}).}
Suppose that $X$ is induced by a $\Z G$-lattice $L$ of type 
$(r,0,0)$. Then 
\[
H^{k}(X/G;\Z)\cong \Z^{\alpha_{k}}\oplus (\Z/p)^{\beta_{k}},
\]
where the coefficients $\alpha_{k}$ and $\beta_{k}$ are given 
as follows: using the same relations as before, 
write the formal power series in $x$ in the form
\[
F_{L}(x)=\left(\frac{1-(\alpha x)^{p}}{1-\alpha x}\right)^{r}
=\sum_{i\ge 0}f_{i}x^{i}+\sum_{i\ge 0}g_{i}\alpha x^{i}
\]
then
\[
\alpha_{k}=\frac{1}{p}\left[\binom{{r(p-1)}}{{k}}+(p-1)(f_{k}-g_{k})
\right]. 
\]
Similarly, $\beta_{k}$ is obtained by writing the formal series 
in $x$ 
\[
T_{L}(x)=\frac{x}{1-x^{2}}\left[
p^{r}x^{2}-x^{2}+1-(1+\alpha x)
\left(\frac{1-(\alpha x)^{p}}{1-\alpha x}\right)^{r}\right]
\]
in the form $T_{L}(x)=\sum_{i\ge 0}\beta_{i}x^{i}+
\sum_{i\ge 0}\gamma_{i}\alpha x^{i}$, where $\alpha$ and 
$\epsilon_{p}$ subject to the same relations as before.
\end{theorem}

\noindent{\bf{Example:}} It is easy to see that in Theorem 
\ref{theoremtype(r,0,0)} the coefficients $\beta_{k}$ are given by
\begin{equation*}
\beta_{k}=\left\{ 
\begin{array}{ccl}
\sum_{j=k}^{r(p-1)}p_{r}(j) & \text{ if } & k 
\text{ is odd and } k>1,\\
0 &  \text{ else.} 
\end{array}%
\right. 
\end{equation*}
Here, as defined above, $p_{r}(j)$ is the number of all
possible sequences of integers $l_{1},...,l_{r}$ such that
$0\le l_{i}\le p-1$ and $l_{1}+\cdots +l_{r}=j$. For example 
when $p=2$, if $X$ is induced by a $\Z G$-lattice of type $(r,0,0)$ 
then 
\begin{equation*}
H^{k}(X/G;\Z)\cong \left\{ 
\begin{array}{ccl}
\Z^{\binom{{r}}{{k}}} &  \text{ if } & k \text{ is even},
\ 0\le k\le r,\\ (\Z/2)^{\binom{{r}}{{k}}+\binom{{r}}{{k+1}}+
\cdots+\binom{{r}}{r}}& \text{ if } & k \text{ is odd},
\ 1<k\le r,\\ 0& \text{ else.} 
\end{array}%
\right. 
\end{equation*}

\medskip

\noindent{\bf{Remark:}} Independently and by different methods, Davis 
and L{\"u}ck in \cite{DavisLuck} obtained the same answers for the 
cohomology groups of the form $H^{*}(X/G;\Z)$, where $X$ is the 
$G$-space induced by a $\Z G$-lattice of type $(r,0,0)$.

\section{Modules of type $(0,s,0)$.}\label{(0,s,0)}

In this section we compute the cohomology groups of the form 
$H^{*}(X/G;\Z)$, where now $X$ is induced by a module $L$ of type 
$(0,s,0)$. We begin by considering the particular case 
$L=(\Z G)^{s}$ and later extend the computations for a general module. 

\begin{lemma}\label{Groupcohomology}
Suppose that $L=(\Z G)^{s}$. Consider the formal power series in $x$
\begin{align*}
F(x)&=(1+\epsilon_{p} x^{p} )^{s}
=\sum_{i\ge 0}f_{i}x^{i}+\sum_{i\ge 0}g_{i}\alpha x^{i},\\
O(x)&=(1+x)^{s}=\sum_{i\ge 0}v_{i}x^{i}
\end{align*}
subject to the same relations as before. Then for $j\ge 1$, 
$H^{j}(X,F)$ is a $\Z G$-lattice of type 
$(g_{j}+v_{j},w_{j}-f_{j},f_{j})$, where
\[
w_{j}=\frac{1}{p}\left[\binom{{sp}}{{j}}-\binom{{s}}{{j}}+
(p-1)(f_{j}-g_{j}-v_{j})\right]. 
\]
\end{lemma}
\Proof
Write $L=P_{1}\oplus \cdots \oplus P_{s}$, where $P_{j}=\Z G$ for 
$1\le j\le s$ and let $i:F\to X$ be the inclusion map. Then, 
up to homeomorphism, $i$ can be identified with 
the diagonal inclusion of $(\S^{1})^{s}$ into $((\S^{1})^{p})^{s}$. 
This shows that 
\[
i^{*}:H^{j}(X)\to H^{j}(F)
\]
is surjective $j\ge 0$. Because of this, the
long exact sequence in cohomology associated
to the pair $(X,F)$ reduces to the different short exact sequences
\[
0\to H^{j}(X,F)\to H^{j}(X)\to H^{j}(F)\to 0.
\]
In particular, $H^{j}(X,F)$ is a free abelian group (ignoring 
the $G$-action) and there is a long exact sequence
\begin{equation}\label{lesZG}
\cdots \to H^{k}(G,H^{j}(X,F))\to H^{k}(G,H^{j}(X))
\stackrel{i^{*}}{\rightarrow} H^{k}(G,H^{j}(F))
\stackrel{\delta}{\rightarrow}\cdots.
\end{equation}
We claim that the map $i^{*}:H^{k}(G,H^{j}(X))\to H^{k}(G,H^{j}(F))$ 
is trivial when $k>0$. To see this note that
\[
H^{j}(X)\cong \wedge^{j} H^{1}(X)\cong
\wedge^{j} L^{*}\cong \wedge^{j} (P_{1}^{*}\oplus \cdots
\oplus P_{s}^{*})
\]
and by the Kunneth theorem there is a commutative diagram
\[
\begin{CD}
H^{k}(G,H^{j}(X)) @>i^{*}>> H^{k}(G,H^{j}(F))\\
@V{\cong}VV @VV{\cong}V\\
\bop_{n_{1}+\cdots n_{s}=j}H^{k}(G,\wedge^{n_{1}}P_{1}^{*}
\otimes\cdots\otimes\wedge^{n_{s}}P^{*}_{s})  @>
\epsilon_{n_{1},...,n_{s}}>>\bop_{n_{1}+\cdots n_{s}=j}H^{k}(G,
\wedge^{n_{1}}\Z\otimes\cdots\otimes\wedge^{n_{s}}\Z),
\end{CD}
\]
where $\epsilon_{n_{1},...,n_{s}}$ is induced by the inclusion
map $i:F\to X$. Notice that the $G$-module $\wedge^{n_{1}}\Z
\otimes\cdots\otimes\wedge^{n_{s}}\Z$ is trivial unless $0\le 
n_{q}\le 1$ for all $q$. Suppose then that $0\le n_{q}\le 1$ for 
all $q$. Since $n_{1}+\cdots+n_{s}=j\ge 1$ it follows that $n_{q}=1$ 
for some $q$ and therefore $\wedge^{n_{1}}P_{1}^{*}\otimes
\cdots\otimes\wedge^{n_{s}}P_{s}^{*}$ has trivial cohomology.
This shows that in any case the map $i^{*}:H^{k}(G,H^{j}(X))\to 
H^{k}(G,H^{j}(F))$ is trivial for $k>0$ and (\ref{lesZG}) reduces 
to the short exact sequence
\begin{equation}\label{sestype0s0}
0\to H^{k-1}(G,H^{j}(F))\stackrel{\delta}{\rightarrow}
H^{k}(G,H^{j}(X,F))
\stackrel{f^{*}}{\rightarrow} H^{k}(G,H^{j}(X))\to 0
\end{equation}
for $k>1$. This sequence splits for $k>1$ as it is a short 
exact sequence of vector spaces over $\F_{p}$. The lemma follows 
using the fact that $H^{j}(X)$ is a $\Z G$-lattice of type $(g_{j},
h_{j}-f_{j},f_{j})$ where 
\[
F(x)=(1+\epsilon_{p}x^{p})^{s}=\sum_{i\ge 0}f_{i}x^{i}+
\sum_{i\ge 0}g_{i}\alpha x^{i}
\]
and the fact that $H^{j}(F)$ is of type $(0,0,v_{j})$ with 
$(1+x)^{s}=\sum_{j\ge 0}v_{j}x^{j}$.
\qed

\medskip

The previous lemma can be used to determine the $p$-torsion 
of $H^{k}_{G}(X,F;\Z)$ in this case. The following 
theorem is then obtained.

\begin{theorem}\label{thmcaseZG}
Suppose that $X$ is induced by the $\Z G$-lattice $L=(\Z G)^{s}$. 
Then the $p$-torsion subgroup of $H^{k}_{G}(X,F;\Z)$ is a 
$\F_{p}$-vector space of dimension $\lambda_{k}$, where 
$\lambda_{k}$ is obtained by writing the formal series in $x$ 
\[
\bar{T}_{L}(x)=\frac{x}{1-x^{2}}\left[
(1+x)^{s}-(1+\alpha x)(1+\epsilon_{p}x^{p})^{s}\right]
\]
in the form $\bar{T}_{L}(x)=\sum_{i\ge 0}\lambda_{i}x^{i}+
\sum_{i\ge 0}\kappa_{i}\alpha x^{i}$, where $\alpha$ and 
$\epsilon_{p}$ subject to the same relations as above.
\end{theorem}
\Proof
As before write $L=P_{1}\oplus \cdots \oplus P_{s}$, where 
$P_{j}=\Z G$ for $1\le j\le s$. The proof of the theorem follows 
by a careful study of the differentials in the Serre spectral 
sequence
\begin{equation}\label{Serre2}
E_{2}^{i,j}=H^{i}(G,H^{j}(X,F))\Longrightarrow H^{i+j}_{G}(X,F).
\end{equation}
Consider the natural map of pairs $f:X\to (X,F)$. If $j>s$ then
\[
f^{*}:H^{j}(X,F)\to H^{j}(X)
\]
is an isomorphism and the Serre spectral sequence
\[
\tilde{E}_{2}^{i,j}=H^{i}(G,H^{j}(X))\Longrightarrow 
H^{i+j}_{G}(X)
\]
collapses on the $E_{2}$-term. Therefore any nontrivial 
differential in (\ref{Serre2}) lands in $E_{k}^{i,j}$ 
for $1\le j\le s$. By the previous lemma for $j\ge 1$
\[
H^{*}(G,H^{j}(X,F))\cong \Z^{w_{j}-f_{j}}\oplus 
(R^{*})^{f_{j}}\oplus (S^{*})^{g_{j}+v_{j}}
\]
where $f_{j},g_{j},w_{j}$ and $v_{j}$ are as described there. 
Note that in particular $v_{j}=\binom{{s}}{{j}}$ for all $j\ge 1$ 
and $f_{j}=g_{j}=0$ if $p\nmid j$. As in the previous 
section all the differentials are trivial in the summand 
$\Z^{w_{j}-f_{j}}\subset H^{0}(G,H^{j}(X,F))$ as it consists 
of elements in the image of the norm map.
We are going to show that in this case all the
differentials are trivial except for the ones of the form
\begin{equation}\label{statement1}
d_{l(p-1)+1}:E_{l(p-1)+1}^{*,lp}\to E_{l(p-1)+1}^{*,l},
\end{equation}
whenever $1\le l\le s$. 
Moreover, we are going to show that
\begin{align}\label{statement2}
E^{*,lp}_{l(p-1)+1}&\cong \Z^{w_{lp}-f_{lp}}
\oplus (R^{*})^{f_{lp}}\oplus (S^{*})^{g_{lp}+v_{lp}}\\
\label{statement3}
E^{*,l}_{l(p-1)+1}&\cong \Z^{w_{l}}\oplus (S^{*})^{v_{l}},
\end{align}
and the homomorphism $d_{l(p-1)+1}$ maps the factor 
$(R^{*})^{f_{lp}}\oplus (S^{*})^{g_{lp}}$ injectively on positive
degrees. Suppose for a moment that this is true. 
Then for all 
$l\ge 1$
\[
q_{E^{*,l}_{l(p-1)+1}}(x)= \frac{v_{l}x}{1-x^{2}}
\]
and since $d_{l(p-1)+1}$ is a homomorphism of graded 
$R^{*}$-modules of degree $l(p-1)+1$ then 
\[
q_{E^{*,l}_{\infty}}(x)=q_{E^{*,l}_{l(p-1)+2}}(x)=
\frac{v_{l}}{1-x^{2}}x-\left(\frac{f_{lp}+g_{lp}x}{1-x^{2}}
\right)x^{l(p-1)+1}.
\]
A similar argument to the one provided below to study the 
differentials in the spectral sequence (\ref{Serre2}) can 
be used to handle the case where the coefficients in the 
sequence (\ref{Serre2}) are $\Q$, $\F_{p}$ and 
$\F_{q}$, where $q$ is a prime different to $p$. This can be 
used to show that there are no extension problems in the 
spectral sequence (\ref{Serre2}). Therefore, the $p$-torsion of 
$H^{k}_{G}(X,F)$ is a $\F_{p}$-vector space of dimension 
$\lambda_{k}$ and 
\begin{align*}
\bar{Q}_{L}(x):&=\sum_{l\ge 0}\lambda_{l}x^{l}
=\sum_{l\ge 0}x^{l}q_{E^{*,l}_{\infty}}(x)\\
&=\sum_{l\ge 1}\left(\frac{v_{l}}{1-x^{2}}\right)x^{l+1}-
\sum_{l\ge 1}\left(\frac{f_{lp}+g_{lp}x}{1-x^{2}}\right)x^{lp+1}\\
&=\frac{x}{1-x^{2}}\left[(1+x)^{s}-
\sum_{l\ge 0}\left(f_{lp}+g_{lp}x\right)
x^{lp}\right].
\end{align*}
Notice that 
\[
F_{L}(x)=(1+\epsilon_{p}x^{p})^{s}=
\sum_{i\ge 0}(f_{i}+\alpha g_{i})x^{i}=
\sum_{l\ge 0}(f_{lp}+\alpha g_{lp})x^{lp},
\]
and thus 
\[
(1+\alpha x)F_{L}(x)=\sum _{i\ge 0}(f_{i}+xg_{i})x^{i}+
\sum _{i\ge 0}\alpha(xf_{i}+g_{i})x^{i}.
\]
Therefore 
\[
\bar{T}_{L}(x):=\frac{x}{1-x^{2}}\left[
(1+x)^{s}-(1+\alpha x)(1+\epsilon_{p}x^{p})^{s}\right],
\]
can be written in the form $\bar{T}_{L}(x)=\bar{Q}_{L}(x)
+\alpha\bar{R}_{L}(x)$, for some formal power series with 
integer coefficients $\bar{R}_{L}(x)$. This together with 
Corollary \ref{free part} proves the theorem.
\smallskip

Induction on $s$ will be used to prove statements 
(\ref{statement1}), (\ref{statement2}) and (\ref{statement3}). 
When $s=1$ then $v_{1}=1$, $f_{p}=1$ and $g_{p}=0$ for $p>2$.
When $p=2$, $f_{2}=0$ and $g_{2}=1$. In any case, by 
Lemma \ref{Groupcohomology} 
the only possible nontrivial differential is
\[
d_{p}:E_{p}^{*,p}\cong \Z^{w_{p}-f_{p}}\oplus 
(R^{*})^{f_{p}}\oplus
(S^{*})^{g_{p}} 
\to E_{p}^{*,1}\cong \Z^{w_{1}}\oplus S^{*}.
\]
In this case, Lemma \ref{cohopairsZG} implies that $H^{i}(X,F)=0$
if $i>sp$. Using Lemma \ref{lemmaBredon} for such $i$ we have
$H^{i}_{G}(X,F)=0$. Since there are no nontrivial differentials
landing in $E^{*,p}_{p}$ it follows that $d_{p}$ is injective for
sufficiently high degrees and the statements follows as $d_{p}$ 
is a homomorphism of graded $R^{*}$-modules of degree $p$.
Suppose that the statements are true for $s\ge 1$ and suppose  
\[
L=P_{1}\oplus\cdots \oplus P_{s+1},
\]
where $P_{i}=\Z G$ for all $1\le j\le s+1$. For every
$1\le k\le s+1$ let
\[
L(k):=P_{1}\oplus \cdots\oplus \hat{P}_{k}\oplus \cdots \oplus 
P_{s+1},
\]
where as usual $\hat{P}_{k}$ means that the factor $P_{k}$ is not
included. For such $k$ we have natural projection and inclusion 
maps
\[
i_{k}:L(k)\to L \text{ and } \pi_{k}:L\to L(k).
\]
These are homomorphisms of $\Z G$-modules that satisfy
$\pi_{k}\circ i_{k}=1$. Moreover, these maps induce 
$G$-equivariant maps
\[
i_{k}:X_{L(k)}\to X_{L} \text{ and } \pi_{k}:X_{L}\to X_{L(k)}
\]
such that that $\pi_{k}\circ i_{k}=1$. By comparing the Serre 
spectral sequences associated to the pairs $(X_{L(k)},F_{L(k)})$ 
and $(X_{L},F_{L})$ using the maps $\pi_{k}$ and $i_{k}$ for 
$1\le k\le s+1$, it follows by induction that
\[
d_{l(p-l)+1}:E_{l(p-1)+1}^{*,lp}\to E_{l(p-1)+1}^{*,l}
\]
is as claimed in (\ref{statement1}) for $1\le l\le s$. Also, we 
can conclude that if $1\le j\le sp$ and $p\nmid j$, then
\[
d_{r}:E^{*,j}_{r}\to E^{*,j-r+1}_{r}
\]
is trivial for all $r$. Consider now the graded $R^{*}$-modules 
of the form $E_{r}^{*,s+1}$, for $r\ge 2$. There are two 
possibilities depending whether $p$ divides $s+1$ or not. 
Suppose first that $p\nmid (s+1)$. In 
this case $f_{s+1}=g_{s+1}=0$ and $v_{s+1}=1$. Thus by 
Lemma \ref{Groupcohomology} 
\[
E_{2}^{*,s+1}\cong\Z^{w_{s+1}-f_{s+1}}\oplus (R^{*})^{f_{s+1}}
\oplus (S^{*})^{g_{s+1}+v_{s+1}}
\cong \Z^{w_{s+1}}\oplus S^{*}.
\]
Note that $1\le s+1\le sp$ and $p\nmid (s+1)$, then by the previous 
comment it follows that all the differentials starting at 
$E_{r}^{*,s+1}$ are trivial for all $r\ge 2$. By 
Lemma \ref{lemmaBredon} there are no nontrivial permament 
cocycles in total degree bigger than  
$(s+1)p$. Therefore if $i$ is odd with $i>(s+1)p $ then all the 
elements in $E^{i,s+1}_{2}\ne 0$ must be in the image of some 
differential $d_{r}:E^{i-r,s+r}\to E^{i,s+1}_{r}$ with 
$r\ge 2$. Note that if $sp<j<(s+1)p$ then $E_{2}^{*,j}\cong 
\Z^{w_{j}}$, in particular $E_{2}^{*,j}$ is a graded 
$R^{*}$-module concentrated on degree $0$ and consists 
of elements in the image of the norm map. Thus all 
differentials starting at $E_{r}^{*,j}$ are trivial 
for all $r\ge 2$. Because of this and the induction hypothesis, 
the only possible nonzero differential that has 
$E_{r}^{*,s+1}$ as target for some $r$ is
\[
d_{(s+1)(p-l)+1}:E_{(s+1)(p-1)+1}^{*,(s+1)p}\to 
E_{(s+1)(p-1)+1}^{*,s+1}.
\]
This shows that $d_{(s+1)(p-l)+1}$ must be nontrivial,
\[
E_{(s+1)(p-l)+1}^{*,(s+1)p}\cong E_{2}^{*,(s+1)p}
\cong \Z^{w_{(s+1)p}-f_{(s+1)p}}\oplus (R^{*})^{f_{(s+1)p}}
\oplus (S^{*})^{g_{(s+1)p}}
\]
and the restriction of $d_{l(p-1)+1}$ to the factor  
$(R^{*})^{f_{(s+1)p}}\oplus (S^{*})^{g_{(s+1)p}}$ is injective 
on positive degrees. This proves the induction hypothesis in 
this case. The case $p\mid (s+1)$ is handled in a similar way.
\qed

\medskip

Suppose now that $L$ is a $\Z G$-lattice of type $(0,s,0)$. 
Then we can find  a $\Z G$-homomorphism 
\[
f:(\Z G)^{s}\to L
\] 
that is an isomorphism after tensoring with $\Z_{(p)}$. By 
Lemma \ref{isocohomologypairs} $f$ induces an isomorphism 
\[
f^{*}:H^{*}_{G}(X_{L},F_{L};\Z_{(p)})\to 
H^{*}_{G}(X_{(\Z G)^{s}},F_{(\Z G)^{s}};\Z_{(p)}).
\]
This shows that the $p$-torsion of $H^{*}_{G}(X_{L},F_{L};\Z)$
is a $\F_{p}$-vector space of dimension $\lambda_{k}$, where 
$\lambda_{k}$ is determined by the formal power series 
$\bar{Q}_{L}(x)=\sum_{j\ge 0}\lambda_{j}x^{j}$. Finally,in 
Lemma \ref{long exact torsion} the long exact sequence 
in cohomology associated to the pair $(X/G,F)$ is studied to 
conclude that $p$-torsion of $H^{k}(X/G;\Z)$ is also 
$\F_{p}$-vector space of dimension $\beta_{k}$ and the 
coefficients $\beta_{k}$ are determined by 
\[
Q_{L}(x):=\sum_{j\ge 0}\beta_{j}x^{j}=\bar{Q}_{L}(x)-
x[(1+x)^{s}-1].
\]
This proves the following theorem.

\begin{theorem}\label{thmcaseZG1}
Suppose that $X$ is induced by a $\Z G$-lattice $L$ of type 
$(0,s,0)$. Then 
\[
H^{k}(X/G;\Z)\cong \Z^{\alpha_{k}}\oplus (\Z/p)^{\beta_{k}},
\]
where the coefficients $\alpha_{k}$ and $\beta_{k}$ are given 
as follows: write the formal power series in $x$
\[
F_{L}(x)=\left(1+\epsilon_{p}x^{p}\right)^{s}=
\sum_{i\ge 0}f_{i}x^{i}+\sum_{i\ge 0}g_{i}\alpha x^{i}
\]
where $\epsilon_{2}=\alpha$ and $\epsilon_{p}=1$ for $p>2$. Then
\[
\alpha_{k}=\frac{1}{p}\left[\binom{{sp}}{{k}}+(p-1)(f_{k}-g_{k})
\right]. 
\]
Similarly, $\beta_{k}$ is obtained by writing the formal series 
in $x$ 
\[
T_{L}(x)=\frac{x}{1-x^{2}}\left[
x^{2}(1+x)^{s}-x^{2}+1-(1+\alpha x)(1+\epsilon_{p} x^{p} )^{s}
\right]
\]
in the form $T_{L}(x)=\sum_{i\ge 0}\beta_{i}x^{i}+
\sum_{i\ge 0}\gamma_{i}\alpha x^{i}$, where $\alpha$ and 
$\epsilon_{p}$ subject to the same relations as above.
\end{theorem}

\section{The general case}\label{general case}

The information collected in the previous two sections is
now assembled to compute the cohomology groups $H^{*}(X/G;\Z)$,
where the $G$-space $X$ is induced by a general $\Z G$-lattice $L$.

\medskip

We start by computing the following special case.

\begin{theorem}\label{equicohtype(r,s,0)}
Suppose that $X$ is the $G$-space associated to a $\Z G$-lattice 
$L$ of type $(r,s,0)$. Then the $p$--torsion subgroup of 
$H^{k}_{G}(X,F;\Z)$ is a $\F_{p}$-vector space of dimension 
$\lambda_{k}$. The coefficients $\lambda_{k}$ are obtained by 
writing the formal series 
in $x$ 
\[
\bar{T}_{L}(x)=\frac{x}{1-x^{2}}\left[
(1+x)^{s}(p^{r}x^{2}-x^{2}+1)-
(1+\alpha x)(1+\epsilon_{p} x^{p} )^{s}
\left(\frac{1-(\alpha x)^{p}}{1-\alpha x}\right)^{r}\right]
\]
in the form $\bar{T}_{L}(x)=\sum_{i\ge 0}\lambda_{i}x^{i}+
\sum_{i\ge 0}\kappa_{i}\alpha x^{i}$, where $\alpha$ and 
$\epsilon_{p}$ subject to the same relations as before.
\end{theorem}
\Proof
For simplicity we are going to consider the case $p>2$. 
The case $p=2$ is handled in a similar way. It is enough 
to prove the theorem for a $\Z G$-lattice of the form 
$L=L_{1}\oplus L_{2}$ where $L_{1}= \oplus_{r} A_{i}$ and 
$L_{2}=\oplus_{s} P_{j}$ with $A_{i}$ an indecomposable 
module of rank $p-1$ as before and $P_{j}=\Z G$ for all 
$1\le j\le s$. This reduction is possible by Lemma 
\ref{isocohomologypairs} and the fact that given any 
$\Z G$-lattice $M$ of type $(r,s,0)$ then we can find a 
homomorphism of $G$-modules $f:L\to M$ that is an isomorphism 
after tensoring with $\Z_{(p)}$. In this case, by 
(\ref{additiveprop}) $X=X_{L_{1}}\times X_{L_{2}}$ and 
$F=F_{L_{1}}\times F_{L_{2}}$. Also $F_{L_{1}}$ is a 
finite set with $p^{r}$ points, $F_{L_{2}}\cong (\S^{1})^{s}$ 
and the inclusion $i_{2}:F_{L_{2}}\to X_{L_{2}}$ can be \
identified, up to homeomorphism, with the diagonal inclusion 
$\Delta:(\S^{1})^{s}\to ((\S^{1})^{p})^{s}$.
Consider the Serre spectral sequence
\begin{equation}\label{Serre3}
E_{2}^{i,j}:=H^{i}(G,H^{j}(X,F;\Z))
\Longrightarrow H^{i+j}_{G}(X,F;\Z).
\end{equation}
We are going to study this spectral sequence by determining 
explicitly it's $E_{2}$-term and the nontrivial differentials. 
We start investigating the $\Z G$-module $H^{j}(X,F)$ 
for $j\ge 1$. When $j=1$ it is easy to see that 
\[
H^{1}(X,F)\cong H^{1}(X_{L_{1}},F_{L_{1}})\oplus 
H^{1}(X_{L_{2}},F_{L_{2}}),
\]
in particular, by the work done in Sections \ref{(r,0,0)} and 
\ref{(0,s,0)} it follows that $H^{1}(X,F)$ is a $\Z G$-lattice of 
type $(s,r,p^{r}-r-1)$. Suppose now that $j\ge 2$ and consider 
the long exact sequence
\[
\cdots\to H^{j-1}(F)\stackrel{\delta}{\rightarrow} H^{j}(X,F)\to
H^{j}(X)\stackrel{i^{*}}{\rightarrow}H^{j}(F)\to \cdots
\]
induced by the pair $(X,F)$. Note that $H^{j}(F)\cong 
\Z^{p^{r}\binom{{s}}{{j}}}$ and $\text{Im}(i^{*})\cong 
\Z^{\binom{{s}}{{j}}}$, therefore the previous sequence reduces 
to the different short exact sequences
\begin{equation}\label{sss(r,s,0)1}
0\to \Z^{(p^{r}-1)\binom{{s}}{{j-1}}}\to H^{j}(X,F)\to M_{j} \to 0,
\end{equation}
where $M_{j}$ fits into the short exact sequence of $\Z G$-modules
\begin{equation}\label{sss(r,s,0)2}
0\to M_{j}\to H^{j}(X)\stackrel{i^{*}}{\rightarrow}
\Z^{\binom{{s}}{{j}}}\to 0.
\end{equation}
Since $M_{j}$ is a subgroup of a free abelian group, then $M_{j}$ 
is a $\Z G$-lattice. To determine it's type recall that $H^{j}(X)$ 
is of type $(g_{j},h_{j}-f_{j},f_{j})$, where
\[
F_{L}(x)=(1+\epsilon_{p} x^{p} )^{s}
\left(\frac{1-(\alpha x)^{p}}{1-\alpha x}\right)^{r}
=\sum_{i\ge 0}f_{i}x^{i}+\sum_{i\ge 0}g_{i}\alpha x^{i}
\]
and $h_{j}$ is determined by $f_{j}$, $g_{j}$ and the rank 
of $H^{j}(X)$. Explicitly,
\[
h_{j}=\frac{1}{p}\left[\binom{{n}}{{j}}+(p-1)(f_{j}-g_{j})
\right]. 
\]
Using the same method that was used in Lemma \ref{Groupcohomology}, 
it can be proved that 
\[
i^{*}:H^{k}(G,H^{j}(X))\to H^{k}(G,\Z^{\binom{{s}}{{j}}})
\]
is the trivial map for $k>0$. Therefore, the long exact sequence 
associated in group cohomology to the short exact sequence 
(\ref{sss(r,s,0)2}) reduces to the short exact sequences 
\[
0\to H^{k-1}(G,\Z^{\binom{{s}}{{j}}})\to H^{k}(G,M_{j})\to 
H^{k}(G,H^{j}(X))\to 0,
\]
for $k>1$. This shows that $M_{j}$ is of type 
$(g_{j}+\binom{{s}}{{j}},m_{j}-f_{j},f_{j})$ and 
$m_{j}$ is determined by $f_{j}$, $g_{j}$ and the rank of $M_{j}$ 
as above. On the other hand, the short exact sequence 
(\ref{sss(r,s,0)1}) yields a long exact sequence
\begin{equation}\label{caseIGZG}
\cdots\to H^{k}(G,\Z^{(p^{r}-1)\binom{{s}}{{j-1}}})\to 
H^{k}(G,H^{j}(X,F))\to H^{k}(G, M_{j})
\stackrel{\partial}{\rightarrow}\cdots.
\end{equation}
When $j\ge 2$  the connecting homomorphism
\[
\partial:H^{k}(G, M_{j})\to 
H^{k+1}(G,\Z^{(p^{r}-1)\binom{{s}}{{j-1}}})
\]
is trivial for $k>0$. This can be seen by comparing the 
long exact sequence (\ref{caseIGZG}) with that associated to 
the pair $(X_{L_{1}}\times F_{L_{2}},F_{L_{1}}\times F_{L_{2}})$. 
Therefore (\ref{caseIGZG}) reduces to the short exact sequence 
\[
0\to H^{k}(G,\Z^{(p^{r}-1)\binom{{s}}{{j-1}}})\to
H^{k}(G,H^{j}(X,F))\to H^{k}(G, M_{j})\to 0
\]
for $k>0$. This sequence splits as it is a short exact sequence 
of $\F_{p}$-vector spaces. This shows that $H^{j}(X,F)$ is a 
$\Z G$-lattice of type
\[
\left(g_{j}+v_{j},w_{j}-f_{j}-u_{j},f_{j}+u_{j}\right),
\]
where $u_{j}=(p^{r}-1)\binom{{s}}{{j-1}}$ and $v_{j}=\binom{{s}}{{j}}$ 
and the coefficient $w_{j}$ is determined in the same way as $h_{j}$.
In particular, there is an isomorphism of graded $R^{*}$-modules 
\begin{align*}
E_{2}^{*,1}&\cong \Z^{r}\oplus (R^{*})^{p^{r}-r-1}\oplus (S^{*})^{s}\\
&\cong \Z^{r}\oplus (R^{*})^{u_{1}-r}\oplus (S^{*})^{v_{1}}
\end{align*}
and for $j\ge 2$
\[
E_{2}^{*,j}\cong \Z^{w_{j}-f_{j}-u_{j}}\oplus (R^{*})^{f_{j}+u_{j}}
\oplus (S^{*})^{g_{j}+v_{j}}.
\]
This completely  describes the $E_{2}$-term of the spectral sequence
(\ref{Serre3}). The differentials in this spectral sequence
can be determined explicitly in a similar fashion as it was done in
Theorems \ref{theoremtype(r,0,0)} and \ref{thmcaseZG}. As a general 
rule, on positive degrees the terms of the form 
$(R^{*})^{f_{j}}\oplus (S^{*})^{g_{j}}$ are the source of a 
nontrivial differential, hence these terms do not survive to 
the $E_{\infty}$-term.  Also the terms of the form  
$(R^{*})^{u_{j}}\oplus (S^{*})^{v_{j}}$ are the target of nontrivial 
differentials. More precisely, write 
\[
\left(\frac{1-(\alpha x)^{p}}{1-\alpha x}\right)^{r}
=\sum_{i\ge0}y_{i}x^{i}+\sum_{i\ge0}z_{i}\alpha x^{i}.
\]
Then an inductive argument on $s$ while keeping $r$ fixed similar 
to the one provided in Theorem \ref{thmcaseZG} shows that
differentials in this spectral sequence are given as follows:
\begin{itemize}
\item
Suppose that $1\le l\le s$, then just as in 
Theorem \ref{thmcaseZG}
\[
d_{l(p-1)+1}:E^{*,lp}_{l(p-1)+1}\to  E^{*,l}_{l(p-1)+1}
\]
is such that 
\[
q_{\text{Im}d_{l(p-1)+1}}(x)=\frac{\binom{{s}}{{l}}}{1-x^{2}}
x^{l(p-1)+1}.
\]

\item Suppose that $1\le l\le s+1$ and $1\le k\le rp$  
are such that $(l-1)(p-1)+k\ge 2$, then
\[
d_{(l-1)(p-1)+k}:E^{*,(l-1)p+k}_{(l-1)(p-1)+k}\to 
E^{*,l}_{(l-1)(p-1)+k}
\]
is such that 
\[
q_{\text{Im}d_{(l-1)(p-1)+k}}(x)=\binom{{s}}{{l-1}}
\left(\frac{y_{k}+z_{k}x}{1-x^{2}}\right)x^{(l-1)(p-1)+k}.
\]

\item All the other differentials are trivial.
\end{itemize}

Let $h(x)=\sum_{k\ge 1}(y_{k}+z_{k}x)x^{k}$. Then 
\begin{align*}
q_{E_{\infty}^{*,1}}(x)&=\frac{1}{1-x^{2}}\left[ 
(sx-sx^{p})+(p^{r}-r-1)x^{2}-\sum_{k\ge 2}(y_{k}+z_{k}x)x^{k}
\right]\\
&=\frac{1}{1-x^{2}}\left[ 
s(x-x^{p})+(p^{r}-1)x^{2}-\sum_{k\ge 1}(y_{k}+z_{k}x)x^{k}
\right]\\
&=\frac{1}{1-x^{2}}\left[ 
s(x-x^{p})+(p^{r}-1)x^{2}-h(x)\right].
\end{align*}
Also, for $j\ge 2$
\begin{align*}
q_{E_{\infty}^{*,j}}(x)&=\frac{1}{1-x^{2}}\left[
v_{j}x+u_{j}x^{2}-\binom{{s}}{{l}}x^{l(p-1)+1}\right]\\ 
&-\frac{1}{1-x^{2}}\binom{{s}}{{j-1}}\left[\sum_{k\ge 1}
(y_{k}+z_{k}x)x^{(j-1)(p-1)+k}\right]\\
&=\frac{1}{1-x^{2}}\left[ 
\binom{{s}}{{j}}\left(x-x^{j(p-1)+1}\right)
+\binom{{s}}{{j-1}}\left((p^{r}-1)x^{2}-h(x)
x^{(j-1)(p-1)}\right)\right].
\end{align*}
In the spectral sequence (\ref{Serre3}) there are no extensions 
problems. This can be seen by studying in the same way the 
sequence (\ref{Serre3}) with coefficients in $\Q$, $\F_{p}$ and 
$\F_{q}$, for a prime $q$ different from $p$. Therefore 
the $p$-torsion of $H_{G}^{k}(X,F;\Z)$ is a $\F_{p}$-vector 
space of dimension $\lambda_{p}$ and these coefficients 
are determined by the formal power series 
\begin{align*}
\bar{Q}_{L}(x)&:=\sum_{j\ge 0}\lambda_{j}x^{j}
=\sum_{j\ge 1}x^{j}q_{E_{\infty}^{*,j}}(x)\\
&=\frac{1}{1-x^{2}}\left[ \sum_{j\ge 1}
\binom{{s}}{{j}}\left(x^{j+1}-x^{jp+1}\right)\right]\\
&+\frac{1}{1-x^{2}}\left[(p^{r}-1)x^{2}\left(\sum_{j\ge 1}
\binom{{s}}{{j-1}}x^{j}\right)-h(x)\left(\sum_{j\ge 1}
\binom{{s}}{{j-1}}x^{(j-1)p+1}\right)
\right]\\
&=\frac{x}{1-x^{2}}\left[((p^{r}-1)x^{2}+1)(1+x)^{s}
-(1+x^{p})^{s}(h(x)+1)\right].
\end{align*}
Note that $h(x)+1=\sum_{k\ge 0}(y_{k}+z_{k}x)x^{k}$
and  
\[
\left(\frac{1-(\alpha x)^{p}}{1-\alpha x}\right)^{r}
=\sum_{i\ge0}(y_{i}+\alpha z_{i})x^{i}.
\]
Therefore 
\[
(1+\alpha x)\left(\frac{1-(\alpha x)^{p}}{1-\alpha x}\right)^{r}
=\sum_{k\ge 0}(y_{k}+z_{k}x)x^{k}+\sum_{k\ge 0}
\alpha(y_{k}x+z_{k})x^{k}. 
\]
This shows that 
\[
\bar{T}_{L}(x):=\frac{x}{1-x^{2}}\left[
(1+x)^{s}((p^{r}-1)x^{2}+1)-(1+\alpha x)(1+\epsilon_{p} x^{p} )^{s}
\left(\frac{1-(\alpha x)^{p}}{1-\alpha x}\right)^{r}\right]
\]
can be written in the form $\bar{T}_{L}(x)=\bar{Q}_{L}(x)
+\alpha\bar{R}_{L}(x)$, for some formal power series with 
integer coefficients $\bar{R}_{L}(x)$.
\qed
\medskip

\begin{lemma}\label{long exact torsion}
Let $X$ be the $G$-space induced by a $\Z G$-lattice $L$ of type 
$(r,s,0)$. Then the $p$-torsion subgroup of $H^{k}(X/G;\Z)$ 
is a finite dimensional vector space over $\F_{p}$ 
of dimension $\beta_{k}$. The coefficients $\beta_{k}$, are 
determined by the formal power series 
\[
Q_{L}(x):=\sum_{k\ge 0}\beta_{k}x^{k}=\bar{Q}_{L}(x)-
x[(1+x)^{s}-1].
\]
Here $\bar{Q}_{L}(x):=\sum_{k\ge 0}\lambda_{k}x^{k}$ and 
$\lambda_{k}$ is the dimension of the $p$-torsion subgroup 
of $H^{k}_{G}(X,F;\Z)$ as a $\F_{p}$-vector space.
\end{lemma}
\Proof
Let us consider first the particular case where the $\Z G$--lattice 
$L$ is of the form $L=(I G)^{r}\oplus (\Z G)^{s}$. The previous 
lemma shows that the $p$-torsion subgroup of $H^{k}_{G}(X,F;\Z)$ is 
a $\F_{p}$-vector space of dimension $\lambda_{k}$. The natural map 
\[
\phi:X\times_{G}EG\to X/G.
\]
induces an isomorphism
\begin{equation}\label{isopairs}
\phi^{*}:H^{*}(X/G,F;\Z)\to H^{*}_{G}(X,F;\Z)
\end{equation} 
by \cite[Proposition VII 1.1]{Bredon}. Therefore the same is 
true for $H^{k}(X/G,F;\Z)$. To handle $H^{*}(X/G;\Z)$ consider 
the long exact sequence in cohomology associated to the pair 
$(X/G,F)$ 
\begin{equation}\label{longexactpair}
\cdots\to H^{j-1}(F;\Z)\stackrel{{\delta}}{\rightarrow}
H^{j}(X/G,F;\Z)\to H^{j}(X/G;\Z)\to H^{j}(F;\Z)\to \cdots.
\end{equation}
By Corollary \ref{fixedpointset} 
\[
F\cong \bigsqcup_{p^{r}}(\S^{1})^{s},
\]
in particular $H^{j}(F;\Z)\cong \Z^{p^{r}\binom{{s}}{{j}}}$. 
By comparing the long exact sequence (\ref{longexactpair})
with the exact sequences in cohomology and equivariant cohomology 
associated to the pair $(X,F)$, it can be proved that 
(\ref{longexactpair}) reduces to the following exact sequences 
\[
0\to \Z^{p^{r}-1}\to H^{1}(X/G,F)\to H^{1}(X/G)\to 
\Z^{s}\to 0,
\]
for $j=1$ and 
\[
0\to(\Z/p)^{\binom{{s}}{{j-1}}}\oplus \Z^{(p^{r}-1)
\binom{{s}}{{j-1}}} \stackrel{\delta}{\rightarrow}
H^{j}(X/G,F)\to H^{j}(X/G)\to \Z^{\binom{{s}}{{j}}}\to 0,
\]
for $j\ge 2$. Also it follows that $\Z^{(p^{r}-1)\binom{{s}}{{j-1}}}$
splits off $H^{j}(X/G,F;\Z)$. Therefore we conclude that the 
$p$-torsion subgroup of $H^{k}(X/G;\Z)$ is a finite dimensional 
$\F_{p}$-vector space of dimension $\beta_{k}$, where 
$\beta_{0}=\lambda_{0}=0$, $\beta_{1}=\lambda_{1}$ and 
\[
\beta_{k}=\lambda_{k}-\binom{{s}}{{k-1}}
\]
for $k\ge 2$. This shows that 
\[
Q_{L}(x)=\bar{Q}_{L}(x)-x[(1+x)^{s}-1]
\]
and the lemma is true in this case. Suppose now that 
$L$ is a general $\Z G$-lattice of type $(r,s,0)$. 
Then we can find a $\Z G$-homomorphism 
\[
f:L'\to L
\]
that is an isomorphism after tensoring with $\Z_{(p)}$ and with 
$L'$ of the kind discussed above so the lemma is true for $L'$.
Using Lemma \ref{isocohomologypairs} we see that  $f$ induces 
an isomorphism
\[
f^{*}:H^{*}_{G}(X_{L},F_{L};\Z_{(p)})\to H^{*}_{G}(X_{L'},F_{L'};
\Z_{(p)}).
\]
This in turn shows that $f$ induces an isomorphism 
\[
f^{*}:H^{*}(X_{L}/G,F_{L};\Z_{(p)})\to H^{*}(X_{L'}/G,F_{L'};
\Z_{(p)})
\]
by (\ref{isopairs}). Finally, by comparing the long exact 
sequences in cohomology with $\Z_{(p)}$-coefficients associated 
to the pairs $(X_{L}/G, F_{L})$ and $(X_{L'}/G, F_{L'})$ 
it follows that $f$ induces an isomorphism 
\[
f^{*}:H^{*}(X_{L}/G;\Z_{(p)})\to H^{*}(X_{L'}/G;\Z_{(p)})
\]
and thus using the universal coefficient theorem we see that 
the lemma is also true for $L$.
\qed
\medskip

We are now ready to prove the main theorem. 

\medskip

\noindent{\bf{Proof of Theorem \ref{maintheorem}:}}
Suppose that $L$ is a $\Z G$-lattice of type $(r,s,t)$. Then we can 
write $L\cong L'\oplus \Z^{t}$ where $L'$ is of type $(r,s,0)$. By 
(\ref{additiveprop}) we have 
\[
X_{L}\cong X_{L'}\times (\S^{1})^{t},
\] 
where $G$ acts trivially on the factor $(\S^{1})^{t}$. Then the 
Kunneth theorem gives an isomorphism
\[
H^{*}(X_{L}/G)\cong H^{*}_{G}(X_{L'}/G)\otimes 
H^{*}((\S^{1})^{t};\Z).
\]
Using this isomorphism and Corollary \ref{free part}, Theorem 
\ref{equicohtype(r,s,0)} and Lemma \ref{long exact torsion} the 
theorem is proved.
\qed

\end{document}